\documentclass[12pt]{amsart}

\usepackage{amsthm}
\usepackage{amsmath}
\usepackage{amssymb}
\usepackage{latexsym}
\usepackage{enumerate}

\usepackage{graphicx}

\newtheorem{theoremA}{Theorem}

\renewcommand{\thetheoremName}

\newtheorem{theorem}{Theorem}[section]

\newtheorem{proposition}[theorem]{Proposition}
\newtheorem{corollary}[theorem]{Corollary}

\theoremstyle{definition}
\newtheorem{definition}[theorem]{Definition}
\newtheorem{example}[theorem]{Example}

\newtheorem{remark}[theorem]{Remark}

\numberwithin{equation}{section}

\newcommand{\D}{\operatorname{D}}
\newcommand{\Hess}{\operatorname{Hess}}
\newcommand{\dist}{\operatorname{dist}}
\newcommand{\Vol}{\operatorname{Vol}}

\newcommand{\C}{\operatorname{Cap}}

\newcommand{\erre}{\mathbb{R}}
\newcommand{\Hp}{\mathbb{H}}
\newcommand{\Lmod}{\operatorname{L}}

\begin{document}

\title[Parabolicity and Hyperbolicity]{On the characterization of parabolicity and hyperbolicity
of submanifolds}

%\author[S. Markvorsen]{Steen Markvorsen$^{\#}$}
%    Address of record for the research reported here
%\address{Department of Mathematics, Technical University of Denmark.}
%    Current address
%\curraddr{}
%\email{S.Markvorsen@mat.dtu.dk}
%    \thanks will become a 1st page footnote.
%\thanks{$^{\#}$Work partially supported by
%the Danish Natural Science Research Council and DGI grant
%MTM2004-06015-C02-02.}

\author[A. Esteve]{Antonio Esteve}
\address{I.E.S. Alfonso VIII, Cuenca-Departament de Matem\`{a}tiques, Universitat Jaume I, Castellon,
Spain.} \email{aesteve7@gmail.com}
\thanks{}

\author[V. Palmer]{Vicente Palmer*}
\address{Departament de Matem\`{a}tiques- Institut de Matem\'atiques i Aplicacions de Castell\'o, Universitat Jaume I, Castellon,
Spain.} \email{palmer@mat.uji.es}
\thanks{* Work partially supported by
the Caixa Castell\'{o} Foundation, and DGI grant MTM2007-62344.}

%    General info
%\date{January 1, 1994 and, in revised form, June 22, 1994.}

%\dedicatory{This paper is dedicated to [[[]]].}


\maketitle
\bibliographystyle{acm}
\begin{abstract}
We give a set of sufficient and necessary conditions for parabolicity and hyperbolicity of a submanifold with controlled mean curvature in a Riemannian manifold with a pole and with sectional curvatures bounded from above or from below.
\end{abstract}

%%%%%%%%%%%%%%%%%%%%%%%%%%%%%%%%%%%%%%%%%%%%%%%%%%%%%%%%%%%%%%%%%%%%%%%%
%
%       SECTION 1 Introduction
%
%%%%%%%%%%%%%%%%%%%%%%%%%%%%%%%%%%%%%%%%%%%%%%%%%%%%%%%%%%%%%%%%%%%%%%%%
\section{Introduction} \label{secIntro}

To find a geometric description for the parabolicity, (or hyperbolicity) of a Riemannian manifold is a question which lies in a central position inside the function theory on Riemannian manifolds, as we can see in the surveys \cite{Li} and \cite{Gri}. This description can be given as a characterization, or as a {\em sufficient} or a {\em necessary} condition. The geometry involved encompasses concepts as the {\em volume growth} of the manifold, or bounds on its Ricci or sectional curvature, (see \cite{A}, \cite{Mi}, \cite{Gri}, \cite {I1}, \cite {I2}, \cite{ChY}, \cite{Va} or, more recently, \cite{HK}).

In $1935$, L.V. Ahlfors proved in \cite{A} that a rotationally symmetric surface $M^2$ is parabolic if and only if the integral $\int_0^\infty \frac{1}{\text{vol}(S(r))}$ is divergent, being $S(r)$ the geodesic circle of radius $r$ in $M^2$. Based on this result, J. Milnor obtained in \cite{Mi} a decision criterion for the parabolicity/hyperbolicity of a complete rotationally symmetric surface which involves its Gaussian curvature. In \cite{Do}, P. G. Doyle showed how to extend these criterion to complete surfaces having a global geodesic polar coordinate system, (namely, having a pole).

 Ahlfors' result has been generalized by several authors, (see \cite{Gri}), to rotationally symmetric spaces with dimension bigger than two, (the so-called {\em model spaces} which will be presented in Subsection \S 2.2), so we have the following theorem:
 \begin{theoremA}[\cite{A}, \cite{Gri}]\label{thmAhlfors}
Let $M^n_w$ be a complete and non compact model space. Then $M^n_w$ is parabolic, (resp. hyperbolic) if and only if
$$\int_{\rho}^\infty \frac{dr}{w^{n-1}(r)}  =\infty \quad (resp. <\infty)$$
\noindent where the volume of the geodesic spheres of $M^n_w$ is given by $\text{vol}(S^w(r))=w^{n-1}(r)$.
\end{theoremA}

 Finally, K. Ichihara proved in \cite{I1}  that a complete, connected and locally compact $n$-Riemannian manifold is parabolic  if its Ricci curvatures are bounded from below by the corresponding curvatures of a model space which satisfies the Ahlfors' integral divergence condition, and it is hyperbolic provided its sectional curvatures are bounded from above by the corresponding curvatures of a model space which satisfies the Ahlfors' integral convergence condition.

In this paper it is considered a submanifold $S^m$ properly immersed in an ambient manifold $N^n$ which has at least one pole and has its radial sectional curvatures, (namely, the sectional curvatures of the planes containing the radial directions from the pole), bounded from above or from below. 

Then, and continuing the programme started with the papers \cite{MP1}, \cite{MP2} and \cite{MP3}, we are going to stablish a set of sufficient conditions for parabolicity and hyperbolicity of submanifolds, (Theorems \ref{thmMain1} and \ref{thmMain2}). These results encompasses partially  the results in \cite{MP2} and \cite{MP3}, and the techniques used to obtain it are based, as in those papers, in the Hessian and Laplacian comparison theory of restricted distance function, which involves bounds on the mean curvature of the submanifold. 

As a consequence of these results, and using the logical interplay among them and the definitions of hyperbolicity and parabolicity, we have obtained two corollaries, (Corollaries \ref{Cor1} and \ref{Cor2}) , with {\em necessary} conditions for these properties. All these results together pretend to approach to a geometric characterization of parabolicity and hyperbolicity for submanifolds in an ambient manifold with bounded (above or below) sectional curvatures, in the style of Theorem \ref{thmAhlfors}.

The way to prove Theorems \ref{thmMain1} and \ref{thmMain2}, (which are the main results of this paper), consist in the application of the Kelvin-Nevanlinna-Royden Criteria, (see \cite[Theorem 5.1]{Gri}), showing the existence of a compact set in the submanifold with positive capacity, (hyperbolicity), or a precompact set with zero capacity, (parabolicity). This method, (which encompasses the use of the distance function from the pole, restricted to the submanifold), is inspired in the {\em Rayleigh's short-cut method} from the classical theory of electricity, used by J. Milnor in \cite{Mi} and by P. G. Doyle in \cite{Do}. 

On the other hand, it was proved in \cite{MP1} that minimal submanifolds of Cartan-Hadamard manifolds are hyperbolic.  We must remark here that in this result it is excluded the case of minimal surfaces in $\erre^3$, (and in $\erre^n$ in general): for example, while the catenoid is parabolic, the doubly periodic Scherk's surface or the triply periodic Swcharz $\mathcal {P}$-surface are hyperbolic. However, the minimal surfaces of the hyperbolic $3$-space are hyperbolic. 
  
In order to explain this particular behaviour,
it was introduced in the paper \cite{MP3}, (see too \cite{HMP}), some control on the
'radiality' of the submanifold. This 'radiality' means the following, assuming for the sole purpose of this explanation, (the proof of our results is independent of the situation of the pole), that the pole $o$ of the ambient manifold lies in the submanifold $S$:
 when the submanifold $S$ is totally geodesic, then $\nabla^N r=\nabla^S r$ in all points,
and, hence, $\Vert \nabla^S r\Vert =1$. On the other hand, and given the starting
point $o \in S$, from which we are measuring the distance $r$, we know that
$\nabla^N r(o)=\nabla^S r(o)$, so $\Vert \nabla^S r(o)\Vert =1$.
Therefore, the difference  $1 - \Vert \nabla^S r\Vert$ quantifies the radial {\em detour}
of the submanifold with respect the ambient manifold as seen
from the pole $o$. To control this detour locally, we apply the following

\begin{definition}

 We say that the submanifold $S$ satisfies a {\it radial tangency
 condition} at $o\in N$ when we have a smooth positive function,
 $$g: S \mapsto \erre_{+} \,\, ,$$ so that
\begin{equation}\label{tangency}
\mathcal{T}(x) \, = \, \Vert \nabla^S r(x)\Vert
\geq g(r(x)) \, > \, 0  \quad {\textrm{for all}}
\quad x \in S \,\, .
\end{equation}
\end{definition}
\medskip

In Corollary 2.2 of \cite{MP3} was proved that a two-dimensional surface $S^2$ in the Euclidean space with the radial component of its mean curvature $H_S$ bounded from below by $0$ is parabolic if the lower bound for its radial tangency $\mathcal{T}(x)$ is a radial function $g(r)$ which is close to $1$ at infinity. 

By contrast, (as it was pointed out there), the Scherk's doubly periodic minimal surface is a hyperbolic surface in $\erre^3$,  such that its radial tangency, (from any fixed point $o$ in the $(x,y)$-plane), is \lq\lq mostly" close to $1$ at infinity, except for the points in the $(x,y)$-plane itself, where the tangency function is close to $0$.

We can single out the following Corollary \ref{MainCor}, (a particular case of Corollary \ref{Cor2}), which explains partially the particular behaviour of the Scherk's surface. Previously to the statement of this Corollary, we need the precise definition of the radially
weighted component of mean curvature:

\begin{definition}
The {\em{$o$-radial mean convexity}}
$\,\mathcal{C}(x)$ of $S$ in $N$, is defined in
terms of the inner product of $H_{S}$ with the
$N$-gradient of the distance function $r(x)$ as
follows:
\begin{equation*}
\mathcal{C}(x) \,=\, -\langle \nabla^{N}r(x),
H_{S}(x) \rangle, \quad x \in S,
\end{equation*}
where $H_{S}(x)$ denotes the mean curvature
vector of $S$ in $N$.

Note that the $o$-radial mean convexity of a minimal submanifold $S$ is $\mathcal{C}(x)=0\,\,\forall x \in S$.
\end{definition}

\begin{corollary}\label{MainCor}
Let $S^2$ be a properly immersed submanifold of $\erre^3$, such that its $o$-radial mean convexity is nonnegative, namely, $\mathcal{C}(x)\geq0\,\,\forall x \in S$. If $S^2$ is hyperbolic, then 

(i) either there isn't exist a smooth positive function
 $$g: S \mapsto \erre_{+} \,\, ,$$ so that
$\mathcal{T}(x) \, = \, \Vert \nabla^S r(x)\Vert
\geq g(r(x)) \, > \, 0  \quad {\textrm{for all}}
\quad x \in S $
\medskip

(ii) or, in case $S$ satisfies a {\it radial tangency
 condition} at $o\in P$ with smooth positive function
 $g: S \mapsto \erre_{+} $, then $\int_{\rho}^{\infty} r e^{-\int_{\rho}^r \frac{m}{t g^2(t)}dt} dr
  < \infty$.
\end{corollary}
\begin{remark}
The hypothesis $\mathcal{C}(x)\geq0\,\,\forall x \in S$ it is satisfied too by the {\em convex} surfaces in $\erre^3$, (see \cite{Pa1}).

On the other hand, in \cite{MP3} is showed as an example how the catenoid, a minimal and parabolic surface, satisfies a {\em radial tangency condition} at the origin $\bar 0 \in \erre^3$.

\end{remark}
\begin{example}

As an example of surface where it is easy to see that this result holds, we have Schwarz P-surface $\mathcal{P} \subseteq \erre^3$. This is a triply periodic minimal surface which is hyperbolic. Its unit cell, (constructed by solving the Plateau problem for a square with corners at the vertices of a regular octahedron), can be viewed roughly as a sphere $S^2$ from which it have been removed six spherical caps whose centroids are antipodal in pairs. In the web page \cite{W}, we can see an image of this unit cell, with the surface-generating straight boundary lines.

We are going to see that assertion (i) of Corollary \ref{MainCor} holds for this surface. To do that, we must remark first that all the points in the ambient space $\erre^3$ are poles. Then, if we consider the center of our extrinsic balls as the center of one of these spheres-unit cells, there exist at least eight points on the surface of this unit cell, (the points where three of the generating straight lines intersect), where $\nabla^{\erre^3} r$ is orthogonal to $\mathcal{P}$. Hence, assertion (i) of Corollary \ref{MainCor} is satisfied. 

\end{example}

%%%%%%%%%%%%%%%%%%%%%%%%%%%%%%%%%%%%%%%%%%%%%%%%%%%
%       SUBSECTION  Outline
%%%%%%%%%%%%%%%%%%%%%%%%%%%%%%%%%%%%%%%%%%%%%%%%%%%
\subsection{Outline of the paper} \label{subsecOutline}
We shall present the basic definitions and results which are in the foundations of our developments in Section 2. Section 3 is devoted to the statement of main theorems and its corollaries.  Proofs of main theorems \ref{thmMain1} and \ref{thmMain2} are presented in Sections 4 and 5.

\subsection{Acknowledgements}

We would like to acknowledge professor Steen
Markvorsen their useful comments concerning these results.

%%%%%%%%%%%%%%%%%%%%%%%%%%%%%
%   Section 2: Preliminaires
%%%%%%%%%%%%%%%%%%%%%%%%%

\section{Preliminaires} \label{Prelim}  We assume
throughout the paper that $S^{m}$ is a non-compact, properly immmersed,
Riemannian submanifold of a complete
Riemannian manifold $N^{n}$. Furthermore, we
assume that $N^{n}$
possesses at least one pole. Recall that a pole
is a point $o$ such that the exponential map
$\exp_{o}\colon T_{o}N^{n} \to N^{n}$ is a
diffeomorphism. For every $x \in N^{n}\setminus \{o\}$ we
define $r(x) = \dist_{N}(o, x)$, and this
distance is realized by the length of a unique
geodesic from $o$ to $x$, which is the {\it
radial geodesic from $o$}. We also denote by $r$
the restriction $r\vert_S: S\to \erre_{+} \cup
\{0\}$. This restriction is called the
{\em{extrinsic distance function}} from $o$ in
$S^m$. The gradients of $r$ in $N$ and $S$ are
denoted by $\nabla^N r$ and $\nabla^S r$,
respectively. Let us remark that $\nabla^S r(x)$
is just the tangential component in $S$ of
$\nabla^N r(x)$, for all $x\in S$. Then we have
the following basic relation:
\begin{equation}\label{radiality}
\nabla^N r = \nabla^S r +(\nabla^N r)^\bot ,
\end{equation}
where $(\nabla^N r)^\bot(x)$ is perpendicular to
$T_{x}S$ for all $x\in S$.

%%%%%%%%%%%%%%%%%%%%%%%%%%%%%%%
%  Subsection: Curvature restrictions
%%%%%%%%%%%%%%%%%%%%%%%%%%

\subsection{Curvature restrictions and extrinsic balls}
\label{subsecCurvRestrict} 

\begin{definition}\label{defRadCurv}
Let $o$ be a point in a Riemannian manifold $M$
and let $x \in M\setminus\{ o \}$. The sectional
curvature $K_{M}(\sigma_{x})$ of the two-plane
$\sigma_{x} \in T_{x}M$ is then called an
\textit{$o$-radial sectional curvature} of $M$ at
$x$ if $\sigma_{x}$ contains the tangent vector
to a minimal geodesic from $o$ to $x$. We denote
these curvatures by $K_{o, M}(\sigma_{x})$.
\end{definition}

\begin{definition}
\begin{enumerate}
\item The submanifold $S$ is called {\em radially 0-convex} if and only if $\mathcal{C}(x) \geq 0 \,\forall x \in S$. 
This condition is satisfied by {\em convex} hypersurfaces of real space forms $\mathbb{K}^{m}(b)$ of constant curvature $b$, (see \cite{Pa1}), as well as by all {\em minimal} submanifolds.
\medskip

\item The submanifold $S$ is called {\em radially minimal} if and only if $\mathcal{C}(x) = 0 \,\forall x \in S$. 
This condition is satisfied by all {\em minimal} submanifolds.
\end{enumerate}
\end{definition}

\begin{definition}\label{ExtBall}
Given a connected and complete
$m$-dimen\-sional submanifold $S^m$ in a complete
Riemannian manifold $N^n$ with a pole $o$, we
denote
the {\em{extrinsic metric balls}} of
(sufficiently large) radius $R$ and center $o$ by
$D_R(o)$. They are defined as any connected
component of the intersection
$$
B_{R}(o) \cap S =\{x\in S \colon r(x)< R\},
$$
where $B_{R}(o)$ denotes the open geodesic ball
of radius $R$ centered at the pole $o$ in
$N^{n}$.
Using these extrinsic balls we define the
$o$-centered extrinsic annuli
$$
A_{\rho,R}(o)= D_R(o) \setminus \bar D_{\rho}(o)
$$
in $S^m$ for $\rho < R$, where $D_{R}(o)$ is the component
of $B_{R}(o) \cap S$ containing $D_{\rho}(o)$.
\end{definition}

\begin{remark}\label{theRemk0}
We want to point out that the extrinsic domains $D_R(o)$
are precompact sets, (because the submanifold $S$ is properly immersed), 
with smooth boundary $\partial D_R(o)$.  The assumption on the smoothness of
$\partial D_{R}(o)$ makes no restriction. Indeed, 
the distance function $r$ is smooth in $N^{n}\setminus \{o\}$ 
since $N^{n}$ is assumed to possess a pole $o\in N^{n}$. Hence
the restriction $r\vert_S$ is smooth in $S$ and consequently the
radii $R$ that produce smooth boundaries
$\partial D_{R}(o)$ are dense in $\mathbb{R}$ by
Sard's theorem and the Regular Level Set Theorem.
\end{remark}

Upper and lower bounds on $\mathcal{C}(x)$, and $\mathcal{T}(x)$  together with a
suitable control on the $o$-radial sectional
curvatures of the ambient space will eventually
control the Laplacian of restricted radial
functions on $S$.

%%%%%%%%%%%%%%%%%%
%Subsection: warped products and model spaces
%%%%%%%%%%%%%%%%%%%%

\subsection{Warped products and model spaces}\label{subsecWarp}
Warped products are generalized manifolds of revolution, see e.g.
\cite{O'N}. Let $(B^{k}, g_{B})$ and $(F^{l}, g_{F})$ denote two
Riemannian manifolds and let $w \colon B \to \mathbb{R_{+}}$ be a
positive real function on $B$. We assume throughout that $w$ is
at least $C^{2}$.
We consider the product manifold $M^{k+l}= B \times F $
and denote the projections onto the factors by $\pi\colon
M \to B$ and $\sigma\colon M \to F$, respectively. The metric
$g$ on $M$ is then defined by the following $w$-modified
(warped) product metric
\begin{equation*}
 g = \pi^{*}(g_{B}) + (w \circ \pi)^{2} \sigma^{*}(g_{F}).
\end{equation*}
\begin{definition}
The Riemannian manifold $(M, g) = (B^{k} \times F^{l}, g) $ is
called a \textit{warped product} with \textit{warping function} $w$,
base manifold $B$ and fiber $F$. We write as follows: $M_{w}^{m} =
B^{k} \times_{w} F^{l}$.
\end{definition}

\begin{definition}[See \cite{Gri}, \cite{GreW}]\label{defModel}
A $w-$model $M_{w}^{m}$ is a smooth warped product with base $B^{1}
= [0,\Lambda[ \,\subset \mathbb{R}$ (where $0 < \Lambda
\leq  \infty$), fiber $F^{m-1} = \mathbb{S}^{m-1}_{1}$ (i.e. the unit
$(m-1)$-sphere with standard metric), and warping function $w\colon
[0,\Lambda[ \to \mathbb{R}_{+}\cup \{0\}$, with $w(0) = 0$,
$w'(0) = 1$, and $w(r) > 0$ for all $r >  0$. The point
$o_{w} = \pi^{-1}(0)$, where $\pi$ denotes the projection onto
$B^1$, is called the {\em{center point}} of the model space. If
$\Lambda = \infty$, then $o_{w}$ is a pole of $M_{w}^{m}$.
\end{definition}

\begin{proposition}\label{propSpaceForm}
The simply connected space forms $\mathbb{K}^{m}(b)$ of constant
curvature $b$ are $w-$models with warping functions
\begin{equation*}
w(r) = Q_{b}(r) =\begin{cases} \frac{1}{\sqrt{b}}\sin(\sqrt{b}\, r) &\text{if $b>0$}\\
\phantom{\frac{1}{\sqrt{b}}} r &\text{if $b=0$}\\
\frac{1}{\sqrt{-b}}\sinh(\sqrt{-b}\,r) &\text{if $b<0$}.
\end{cases}
\end{equation*}
Note that for $b > 0$ the function $Q_{b}(r)$ admits a smooth
extension to  $r = \pi/\sqrt{b}$.
\end{proposition}

\begin{proposition}[See \cite{O'N},  \cite{GreW} and \cite{Gri}]\label{propWarpMean}
Let $M_{w}^{m}$ be a $w-$model with warping function $w(r)$ and
center $o_{w}$. The distance sphere of radius $r$ and center $o_{w}$
in $M_{w}^{m}$ is the fiber $\pi^{-1}(r)$. This distance sphere has
the constant mean curvature $\eta_{w}(r)= \frac{w'(r)}{w(r)}$ On the other hand, the
$o_{w}$-radial sectional curvatures of $M_{w}^{m}$ at every $x \in
\pi^{-1}(r)$ (for $r > 0$) are all identical and determined
by
\begin{equation*}
K_{o_{w} , M_{w}}(\sigma_{x}) = -\frac{w''(r)}{w(r)}.
\end{equation*}
\end{proposition}

%%%%%%%%%%%%%%%%%%%%%%%%
%Subsection: Hessian and Laplacian comparison
%%%%%%%%%%%%%%%

\subsection{Hessian and Laplacian comparison analysis}\label{subsecLap}
The 2.nd order analysis of the restricted distance function $r_{|_{P}}$ defined on manifolds with a pole is firstly and foremost governed by the Hessian comparison Theorem A in \cite{GreW}:

\begin{theorem}[See \cite{GreW}, Theorem A]\label{thmGreW}
Let $N=N^{n}$ be a manifold with a pole $o$, let $M=M_{w}^{m}$ denote a
$w-$model with center $o_{w}$, and $m \leq n$. Suppose that every $o$-radial
sectional curvature at $x \in N \setminus \{o\}$ is bounded from above by
the $o_{w}$-radial sectional curvatures in $M_{w}^{m}$ as follows:
\begin{equation*}
K_{o, N}(\sigma_{x}) \geq\, (\leq)\, -\frac{w''(r)}{w(r)}
\end{equation*}
for every radial two-plane $\sigma_{x} \in T_{x}N$ at distance $r =
r(x) = \dist_{N}(o, x)$ from $o$ in $N$. Then the Hessian of the
distance function in $N$ satisfies
\begin{equation}\label{eqHess}
\begin{aligned}
\Hess^{N}(r(x))(X, X) &\leq\,(\geq)\, \Hess^{M}(r(y))(Y, Y)\\ &=
\eta_{w}(r)\left(1 - \langle \nabla^{M}r(y), Y \rangle_{M}^{2}
\right) \\ &= \eta_{w}(r)\left(1 - \langle \nabla^{N}r(x), X
\rangle_{N}^{2} \right)
\end{aligned}
\end{equation}
for every unit vector $X$ in $T_{x}N$ and for every unit vector $Y$
in $T_{y}M$ with $\,r(y) = r(x) = r\,$ and $\, \langle
\nabla^{M}r(y), Y \rangle_{M} = \langle \nabla^{N}r(x), X
\rangle_{N}\,$.
\end{theorem}

\begin{remark}
In \cite[Theorem A, p. 19]{GreW}, the Hessian of
$r_M$ is less or equal to the Hessian of $r_N$
provided that the radial curvatures of $N$ are
bounded from above by the radial curvatures of
$M$ and provided that $\dim M \geq \dim N$. This
latter dimension condition is {\em{not}}
satisfied in our setting. However, since $(M^{m},
g)$ is a $w-$model space it has an
$n-$dimensional $w-$model space companion with
the same radial curvatures and the same Hessian
of radial functions as $(M^{m}, g)$. In effect,
therefore, applying \cite[Theorem A, p. 19]{GreW}
to the high-dimensional comparison space gives
the low-dimensional comparison inequality as
stated.

In other words,  $\Hess^{M_w}(r(y))(Y, Y)$ {\em do not}
depend on the dimension $m$, as we can easily  see by computing it directly,
(see \cite{Pa3}), so the hypothesis on the dimension can be overlooked in the comparison among the Hessians.
\end{remark}

As a consecuence of this result, we have the following Laplacian inequalities:
\begin{proposition} \label{corLapComp}
Let $N^{n}$ be a manifold with a pole $p$, let $M_{w}^{m}$ denote a $w-$model
with center $p_{w}$. 
\medskip

(i) Suppose that every $o$-radial sectional curvature at $x \in
N - \{o\}$ is bounded from below by the $o_{w}$-radial sectional curvatures in
$M_{w}^{m}$ as follows:
\begin{equation}\label{eqKbound}
\mathcal{K}(\sigma(x)) \, = \, K_{o,
N}(\sigma_{x}) \geq  -\frac{w''(r)}{w(r)}
\end{equation}
for every radial two-plane $\sigma_{x} \in T_{x}N$ at distance $r = r(x) =
\dist_{N}(o, x)$ from $o$ in $N$. Then we have for
every smooth function $f(r)$ with $f'(r) \leq
0\,\,\textrm{for all}\,\,\, r$, (respectively
$f'(r) \geq 0\,\,\textrm{for all}\,\,\, r$):
\begin{equation} \label{eqLap1}
\begin{aligned}
\Delta^{S}(f \circ r) \, \geq (\leq) \, &\left(\,
f''(r) - f'(r)\eta_{w}(r) \, \right)
 \Vert \nabla^{S} r \Vert^{2} \\ &+ mf'(r) \left(\, \eta_{w}(r) +
\langle \, \nabla^{N}r, \, H_{S}  \, \rangle  \,
\right)  \quad ,
\end{aligned}
\end{equation}
where $H_{S}$ denotes the mean curvature vector
of $S$ in $N$.
\medskip

(ii)  Suppose that every $o$-radial sectional curvature at $x \in
N - \ o\}$ is bounded from above by the $o_{w}$-radial sectional curvatures in
$M_{w}^{m}$ as follows:
\begin{equation}\label{eqKbound}
\mathcal{K}(\sigma(x)) \, = \, K_{o,
N}(\sigma_{x}) \leq  -\frac{w''(r)}{w(r)}
\end{equation}
for every radial two-plane $\sigma_{x} \in T_{x}N$ at distance $r = r(x) =
\dist_{N}(o, x)$ from $p$ in $N$. Then we have for
every smooth function $f(r)$ with $f'(r) \leq
0\,\,\textrm{for all}\,\,\, r$, (respectively
$f'(r) \geq 0\,\,\textrm{for all}\,\,\, r$):
\begin{equation} \label{eqLap1}
\begin{aligned}
\Delta^{S}(f \circ r) \, \leq (\geq) \, &\left(\,
f''(r) - f'(r)\eta_{w}(r) \, \right)
 \Vert \nabla^{S} r \Vert^{2} \\ &+ mf'(r) \left(\, \eta_{w}(r) +
\langle \, \nabla^{N}r, \, H_{S}  \, \rangle  \,
\right)  \quad ,
\end{aligned}
\end{equation}
\end{proposition}
%%%%%%%%%%%%%%%%%%%%%%%%%%%%%%%
%  Subsection: Capacities of extrinsic annular domains
%%%%%%%%%%%%%%%%%%%%%%%%%%
\subsection{Capacities of extrinsic annular domains}\label{subsecCap}

The proof of theorems \ref{thmMain1} and \ref{thmMain2} is based on the existence
of (lower and upper) bounds for the capacity of some compact subset in the submanifold $S^m$.
This compact subset is an extrinsic ball $\D_{\rho}(o) \subseteq S$. 

In general, the capacity of a compact domain $K$ in a
    precompact open set $\Omega$ of a
Riemannian manifold $M$ can be expressed as the following integral along the
boundary of the compact set $K$ (see e.g. \cite{Gri}):
\begin{equation}\label{eq1.1}
\C(K,\Omega)=\int_{\Omega-K}\Vert \nabla^M v\Vert^2 dV \quad ,
\end{equation}
where the function $v$ is the solution of the
Dirichlet problem in $\Omega-K$
\begin{equation}
\label{eqDirPoEq}
\begin{aligned} \begin{cases}
\Delta^M v &= 0\,\,\,\text{on $\Omega - K$}\\
\phantom{\Delta^M  }v &=0\,\,\,\text{on $\partial K$}  \\
\phantom{\Delta^M }v&= 1\,\,\,\text{on $\partial \Omega$} \quad .
\end{cases}
\end{aligned}
\end{equation}

The capacity of $K$ in the whole manifold $M$ is given by the following limit, given any exhaustion sequence of precompact open
    subsets
    $\{\Sigma_n\}_{n \in \mathbb{N}}$
    covering all of $M$ such that $\Sigma_0 = K$ and $\Sigma_n \subseteq \Sigma_{n+1}$, (see \cite{Gri}):
    \begin{equation}
    \C(K,M)=\lim_{n \to \infty} \C(K, \Sigma_{n})
    \end{equation}

Using the divergence theorem, it is easy to see, (\cite{Gri}), that integral (\ref{eq1.1}) becomes

\begin{equation}\label{eq1.2}
\C(K,\Omega)=\int_{\partial K}\langle \nabla^M v,\nu\rangle dA \quad ,
\end{equation}
where $\nu$ is the unit normal vector field on $\partial K$ which points
{\em{into}} the domain $\Omega-K$.
    
    In our setting, we have a compact set in the submanifold $S^m$, $\bar D_{\rho}(o)$, and an exhaustion of $S$ given by the extrinsic balls $\{D_R(o)\}_{R>0}$ which contains $\bar D_{\rho}(o)$. The  computation of the capacity of these extrinsic annular domains $A_{\rho, R}= D_R(o) \setminus \bar D_{\rho}(o)$ is given by the following considerations, applying equation (\ref{eq1.2}):
   
\begin{equation}
\C(A_{\rho, R}) \,= \,  \int_{\partial D_{\rho}} \langle \nabla^{P}v, n_{\partial D_{\rho}} \rangle_{\partial D_{\rho}} \, d\mu \quad,
\end{equation}
where $v(x)$ is the Laplace potential function
for the extrinsic annulus  $A_{\rho, R}\, = \,
D_{R} - D_{\rho}$, setting $v_{\partial D_{\rho}}\,=\,0$ and
$v_{\partial D_{R}}\,=\,1$ and  $n_{\partial D_{\rho}}$ denotes the unit normal vector field along $\partial D_{\rho}$ pointing {\em{into}} the domain
$A_{\rho, R}$. 

The function $v$ must be nonnegative in the annular domain $A_{\rho, R}$. Otherwise $v$ would have an intrinsic (negative) minimum in $A_{\rho, R}$, and since $v$ is harmonic this is ruled out by the minimum principle.

Now, since $v$ is nonnegative and $v=0$ at the inner boundary, then the inwards directed gradient  $\langle \nabla^{P}v, n_{\partial D_{\rho}} \rangle_{\partial D_{\rho}}$ is also nonnegative.
Since $\partial D_{\rho}$ is a level hypersurface (of value $v=0$) for $v$ in $P$, we have that $n_{\partial D_{\rho}}$ is proportional to $\nabla^{P}v$. It therefore follows that
\begin{equation}
\langle \nabla^{P}v, n_{\partial D_{\rho}} \rangle_{\partial D_{\rho}} \,= \, \Vert \nabla^{P}v(x) \Vert \quad.
\end{equation}

Therefore we have
    \begin{equation}\label{CapEq}
     \C(A_{\rho, R}) = \, \int_{\partial D_{\rho}}\Vert \nabla^{P}v(x) \Vert \, d\nu 
         \end{equation}

%%%%%%%%%%%%%%%%%%%%%%%%%%
%Section 3: Main results
%%%%%%%%%%%%%%%%%%%%%%

\section{Main results}\label{secMain}
We are going to give some previous definitions, in order to  
formulate our main hyperbolicity and parabolicity results. The
proofs are developed through the following
sections.

\begin{definition}
Let $N^n$ be a complete manifold with pole $o \in N$, and let $S^m$ be a properly immersed submanifold in $N$. Given a function $h: S \longrightarrow \erre$ which only depends on the extrisic distance $r$ in $S$, $h(r(x))$ for all  $x \in S$, we say that the function $h(r)$ is balanced from above with respect the warping function $w(r)$ of a model space $M^m_w$ if 
\begin{equation}\label{balance1}
\mathcal{M}(r)=m(\eta_w(r)-h(r)) \geq 0\,\,\forall r
\end{equation}
and that the function $h(r)$ is balanced from below with respect the warping function $w(r)$ if 
\begin{equation}\label{balance2}
\mathcal{M}(r)=m(\eta_w(r)-h(r)) \leq 0\,\,\forall r
\end{equation}
\end{definition}

\begin{remark}
As in the  following parabolicity and hyperbolicity criteria plays a fundamental r\^ole the convergence/divergence of the infinite integrals $\int_{\rho}^{\infty} \Lambda_g(t)\,  dt $ and $\int_{\rho}^{\infty} \Lambda(t)
\,  dt $, we should remark that, if $\mathcal{M}(r)\geq\, 0$ for all $r >0$, then $\int_{\rho}^{\infty} \Lambda_g(t)\,  dt \,\leq \,\int_{\rho}^{\infty} \Lambda(t)
\,  dt $ and, on the other hand,  if $\mathcal{M}(r)\leq\, 0$ for all $r >0$, then $\int_{\rho}^{\infty} \Lambda_g(t)\,  dt \,\geq \,\int_{\rho}^{\infty} \Lambda(t)
\,  dt $
\end{remark}
\begin{definition}\label{Lambda}
Let $N^n$ be a complete manifold with pole $o \in N$, and let $S^m$ be a properly immersed submanifold in $N$. Let us consider too a model space $M^m_w$.

(i) Define $\Lambda(r)$ as the function
\begin{equation*}
\Lambda(r) =
w(r)\exp\left(-\int_{\rho}^{r}
\mathcal{M}(t)\,
dt\right).
\end{equation*}
\medskip

(ii) Assume moreover that $S$ satisfies a {\it radial tangency
 condition} at $o\in N$.
 
 We denote as $\Lambda_g(r)$ the function
\begin{equation*}
\Lambda_g(r) =
w(r)\exp\left(-\int_{\rho}^{r}
\frac{\mathcal{M}(t)}{g^{2}(t)}\,
dt\right).
\end{equation*}
\end{definition}

\begin{theorem}[Parabolicity]\label{thmMain1}
Let $N^n$ be a complete manifold with pole $o$, and suppose that 
\begin{equation} \label{eqKcomp1}
K_{o, N}(\sigma_{x}) \geq
-\frac{w''(r)}{w(r)}
\end{equation}
for all $x$ with $r=r(x)\in [0,\infty)$. 

Let $S^m$ be a complete and properly immersed submanifold with $o$-radial mean convexity $\mathcal{C}(x)$ bounded from below by the radial function $h(r(x))$:
\begin{equation}\label{geConvexity}
\mathcal{C}(x) \geq h(r(x)) \,\,\,\textrm{for
all}\,\,\, x \in S^{m} \,\,\,
{\textrm{with}}\,\,\, r(x) \in [0, \infty).
\end{equation}

Then:
\medskip

(A) Assume that the submanifold $S$ satisfies a {\it radial tangency
 condition} at $o\in N$, (namely, there exists smooth 
 $g: S \mapsto \erre_{+} \,\, ,$ so that $\Vert \nabla^S r(x)\Vert
\geq g(r(x)) \, > \, 0 \, {\textrm{for all}}
\quad x \in S$), and that the function $h(r)$ is balanced from above with respect the warping function $w(r)$, ($\mathcal{M}(r) \geq 0\,\forall r$). 
 
 \medskip
Suppose that 
\begin{equation} \label{eqParabCond}
\int_{\rho}^{\infty} \Lambda_g(t)
\,  dt =  \infty.
\end{equation}
Then $S^{m}$ is parabolic.
\bigskip

(B) Assume that  the function $h(r)$ is balanced from below  with respect the warping function $w(r)$ , ($\mathcal{M}(r) \leq 0\,\forall r$), and suppose that 
\begin{equation} \label{eqParabCond2}
\int_{\rho}^{\infty} \Lambda(t)
\,  dt =  \infty.
\end{equation}
Then $S^{m}$ is parabolic.
\end{theorem}

\begin{remark}
Theorem \ref{thmMain1} (A) has been stated and proved in \cite{MP3}, (see Theorem 9.2), under a more restricitive balance condition.
\end{remark}

\begin{theorem}[Hyperbolicity]\label{thmMain2}
Let $N^n$ be a complete manifold with pole $o$, and suppose that 
\begin{equation} \label{eqKcomp2}
K_{o, N}(\sigma_{x}) \leq
-\frac{w''(r)}{w(r)}
\end{equation}
for all $x$ with $r=r(x)\in [0,\infty)$. Let $S^m$ be a complete and properly immersed submanifold with $o$-radial mean convexity $\mathcal{C}(x)$ bounded from above by the radial function $h(r(x))$:
\begin{equation}\label{leConvexity}
\mathcal{C}(x) \leq h(r(x)) \,\,\,\textrm{for
all}\,\,\, x \in S^{m} \,\,\,
{\textrm{with}}\,\,\, r(x) \in [0, \infty).
\end{equation}

Then
\medskip

(A) Assume that the submanifold $S$ satisfies a {\it radial tangency
 condition} at $o\in P$, and that the function $h(r)$ is balanced from below with respect the warping function $w(r)$, ($\mathcal{M}(r) \leq 0\,\forall r$). 
 
 \medskip
Suppose finally that 
\begin{equation} \label{eqHypCond}
\int_{\rho}^{\infty} \Lambda_g(t)
\,  dt < \infty.
\end{equation}
Then $S^{m}$ is hyperbolic.
\bigskip

(B) Assume that  the function $h(r)$ is balanced from above with respect the warping function $w(r)$, ($\mathcal{M}(r) \geq 0\,\forall r$), and suppose that 
\begin{equation} \label{eqHypCond2}
\int_{\rho}^{\infty} \Lambda(t)
\,  dt <  \infty.
\end{equation}
Then $S^{m}$ is hyperbolic.
\end{theorem}

\begin{remark}
Theorem \ref{thmMain2} (B) has been stated and proved in \cite{MP2}.
If we follow the notation in \cite{MP2}, we have
$$\mathcal{G}(r)=\exp(\int_{\rho}^{r}h(t)\,dt)$$
and it is straightorward to check that
 $\int_{\rho}^{\infty}\frac{\mathcal{G}^{m}(r)}{w^{m-1}(r)}\,dr < \infty$ iff $\int_{\rho}^{\infty} \Lambda(t)
\,  dt <  \infty$.
\end{remark}

\begin{remark}
We have the following examples of a direct application of theorems \ref{thmMain1} and \ref{thmMain2}. Concerning Theorem \ref{thmMain1}, we can see as the cones and the paraboloids, (both convex hypersurfaces in $\erre^3$), are parabolic, (see \cite{MP3}). Concerning Theorem \ref{thmMain2}, we have that  surfaces $P^2$ in $\Hp^3(b)$ with constant mean curvature $H_P \leq \frac{1}{2}\sqrt{-b}$ are hyperbolic, (see too Corollary B in  \cite{MP2}).
\end{remark}

%%%%%%%%%%%%%%%%%%%%%%%%
%Subsection: Corollaries
%%%%%%%%%%%%%%%

\subsection{Corollaries}\label{Corol}

Finally, as corollaries of Theorem \ref{thmMain1} and Theorem \ref{thmMain2} we have the following results.

\begin{corollary}\label{Cor2}
Let $N^n$ be a complete manifold with pole $o$, and let $S^m$ be a properly immersed submanifold
in $N$, both satisfying inequalities (\ref{eqKcomp1}) and (\ref{geConvexity}) in Theorem \ref{thmMain1}.

(A) Let us suppose that $\mathcal{M}(r) \geq 0\,\forall r$.

If $S$ is hyperbolic, then 
\medskip

(A.1) either there isn't exist a smooth positive function
 $$g: S \mapsto \erre_{+} \,\, ,$$ so that
$\mathcal{T}(x) \, = \, \Vert \nabla^S r(x)\Vert
\geq g(r(x)) \, > \, 0  \quad {\textrm{for all}}
\quad x \in S $
\medskip

(A.2) or, in case $S$ satisfies a {\it radial tangency
 condition} at $o\in P$ with smooth positive function
 $g: S \mapsto \erre_{+} $, then $\int_{\rho}^{\infty}\Lambda_g(r) dr
  < \infty$.
  
  \medskip
  
  (B) Let us suppose that $\mathcal{M}(r) \leq 0\,\forall r$.
  
  If $S$ is hyperbolic, then 
  $\int_{\rho}^{\infty}\Lambda(r) dr
  < \infty$.
\end{corollary}
\begin{proof} 
(A) The ambient manifold $N$ and the submanifold $S$ satisfies hypothesis (\ref{eqKcomp1}) and (\ref{geConvexity}) in Theorem \ref{thmMain1}.  Hence, if $S$ is not parabolic, then we have the negation of both sets of assumptions in assertions (A) and (B) in Theorem \ref{thmMain1}. In particular, assertion (A) doesn't holds so, as $\mathcal{M}(r) \geq 0\,\forall r$, then either there is not any smooth positive function
 $$g: S \mapsto \erre_{+} \,\, ,$$ so that
$\mathcal{T}(x) \, = \, \Vert \nabla^S r(x)\Vert
\geq g(r(x)) \, > \, 0  \quad {\textrm{for all}}
\quad x \in S $
or, if  this bounding function for the tangency exists, then $\int_{\rho}^{\infty}\Lambda_g(r) dr
  < \infty$.
  \medskip
  
  (B) In this case, $\mathcal{M}(r) \leq 0\,\forall r$, and, as assertion (B) in Theorem \ref{thmMain1} doesn't holds, we  conclude that $\int_{\rho}^{\infty}\Lambda(r) dr
  < \infty$.
\end{proof}

\begin{remark}
Corollary \ref{MainCor} in the Introduction follows from Corollary \ref{Cor2}, if we consider that the ambient manifold $N$ is the Euclidean $3$-space $\erre^3$, (which implies to consider as a warping function $w(r)=r$), and we have into account that, by hypothesis, $h(r(x))=0$ for all $x \in S$, and hence $\mathcal{C}(x)=\eta_w(r(x))-h(r(x))=\frac{1}{r}>0 $ for all $x \in S$.
\end{remark}

\begin{corollary}\label{Cor1}
Let $N^n$ be a complete manifold with pole $o$, and let $S^m$ be a properly immersed submanifold in $N$, both satisfying inequalities (\ref{eqKcomp2}) and (\ref{leConvexity}) in Theorem \ref{thmMain2}.

(A) Let us suppose that $\mathcal{M}(r) \geq 0\,\forall r$.

If $S$ is parabolic, then 
\begin{equation} \label{eqParabCond3}
\int_{\rho}^{\infty} \Lambda(t)
\,  dt =  \infty.
\end{equation}
\medskip

(B) Let us suppose that $\mathcal{M}(r) \leq 0\,\forall r$.

If $S$ is parabolic, then 
\medskip

(B.1) either there isn't exist a smooth positive function
 $$g: S \mapsto \erre_{+} \,\, ,$$ so that
$\mathcal{T}(x) \, = \, \Vert \nabla^S r(x)\Vert
\geq g(r(x)) \, > \, 0  \quad {\textrm{for all}}
\quad x \in S $
\medskip

(B.2) or, in case $S$ satisfies a {\it radial tangency
 condition} at $o\in P$ with smooth positive function
 $g: S \mapsto \erre_{+} $, then $\int_{\rho}^{\infty}\Lambda_g(r) dr
  =\infty$.

\end{corollary}
\begin{proof} 

(A) The ambient manifold $N$ and the submanifold $S$ satisfies hypothesis (\ref{eqKcomp2}) and (\ref{leConvexity}) in Theorem \ref{thmMain2}. Hence, if $S$ is not hyperbolic, then we have the negation of both sets of assumptions in assertions (A) and (B) in Theorem \ref{thmMain2}. In particular, some of the two assumptions in assertion (B) doesn't holds so, as $\mathcal{M}(r) \geq 0\,\forall r$, we have
\begin{equation} 
\int_{\rho}^{\infty} \Lambda(t)
\,  dt= \infty.
\end{equation}
\medskip
(B) In this case,  $\mathcal{M}(r) \leq 0\,\forall r$, and  with same arguments than before, namely, concluding the negation of assertion (A) in Theorem  \ref{thmMain2}, we have that either there isn't exist a smooth positive function
 $$g: S \mapsto \erre_{+} \,\, ,$$ so that
$\mathcal{T}(x) \, = \, \Vert \nabla^S r(x)\Vert
\geq g(r(x)) \, > \, 0  \quad {\textrm{for all}}
\quad x \in S $
\medskip
or, in case $S$ satisfies a {\it radial tangency
 condition} at $o\in P$ with smooth positive function
 $g: S \mapsto \erre_{+} $, then $\int_{\rho}^{\infty}\Lambda_g(r) dr
  =\infty$.
\end{proof}

\begin{corollary}\label{Cor3}
Let $M^n_w$ be a model space, and let $S^m$ be a properly immersed submanifold of $M^n_w$. Assume that the $o$-radial mean convexity of $S$ in $M^n_w$, is  equal to the radial function $\eta_w(r(x))$, namely
\begin{equation}\label{eqConvexity}
\mathcal{C}(x) =\eta_w(r(x)) \,\,\,\textrm{for
all}\,\,\, x \in S^{m} \,\,\,
{\textrm{with}}\,\,\, r(x) \in [0, \infty).
\end{equation}

Then, $S^m$ is parabolic if and only if $\int_{\rho}^{\infty} w(r) \,  dr = \infty$.
\end{corollary}

%%%%%%%%%%%%%%%%%%%%%
%       Section 4 Proof of Theorem \ref{thmMain1} (i)} and Theorem \ref{thmMain2} (i)}
%%%%%%%%%%%%%%%%%%%%%%%%%%%%%%%%%%%%%%%%%%%%%%%%%%%

%\section{Drifted $2$-capacity of model spaces} \label{secDrift}
\section{Proof of  assertion (A) in Theorem \ref{thmMain1} and Theorem \ref{thmMain2}} \label{proofthmMain1}

%{\bf Proof of Theorem \ref{thmMain1} (i)} and Theorem \ref{thmMain2} (i)\\

As the submanifold $S$ satisfies a radial tangency condition, given the function $g: S \longrightarrow \erre_+$, we define a second order differential operator on functions of one real variable as follows:

\begin{equation}\label{operator}
\Lmod_g \psi(r) = \psi''(r) + \psi'(r)\left(
\frac{\mathcal{M}(r)}{g^{2}(r)} -
\eta_{w}(r)\right).
\end{equation}

and consider the smooth solution $\psi_{\rho,R}(r)$ of the following Dirichlet-Poisson problem
associated to $\Lmod_g$:

\begin{equation}\label{eqDir}
\begin{cases}
\Lmod_g \psi &= 0\,\,\,\text{on $[\rho, R]$}\\
\phantom{L }\psi(\rho) &= 0 \\
\phantom{L }\psi(R)&= 1
\end{cases}
\end{equation}

The explicit solution to the Dirichlet problem
\eqref{eqDir} is given in the following
Proposition which is straightforward,

\begin{proposition}\label{propDirSol}
The solution to the Dirichlet problem (\ref{eqDir}) only depends on
$r$ and is given explicitly - via the function
$\Lambda_g(r)$ introduced in Definition \ref{Lambda} (ii), by:
\begin{equation} \label{eqPsi}
\psi_{\rho,R}(r) =  \frac{\int_{\rho}^{r}
\Lambda_g(t)\,dt}{\int_{\rho}^{R}\Lambda_g(t)\,dt}.
\end{equation}
The corresponding 'drifted' capacity is
\begin{equation} \label{eqModelCap}
\begin{aligned}
\C_{\Lmod_g}(A_{\rho, R}^{w})
&=\int_{\partial D_{\rho}^{w}}\langle\nabla^{M}\psi_{\rho,R},\nu\rangle\,dA\\
&=\Vol(\partial D_{\rho}^{w})\Lambda_g(\rho)\left(\int_{\rho}^{R}
\Lambda_g(t)\,dt\right)^{-1}.
\end{aligned}
\end{equation}
\end{proposition}

At this point the proof of these two results splits, in the following way:
\medskip

{\bf Assertion (A) in Theorem \ref{thmMain1}}.
\medskip
\\
Concerning the proof of assertion (A) in Theorem \ref{thmMain1} it is easy to see, using equation (\ref{eqPsi}) and the balance condition (\ref{balance1}) that the solution $\psi_{\rho,R}$ of the problem (\ref{eqDir}) satisifies:
\begin{equation}\label{parenth}
\begin{aligned}
\psi'_{\rho,R}(r)& \geq 0\\
\psi''_{\rho,R}(r)& - \psi'_{\rho,R}(r)\eta_{w}(r) =- \psi'_{\rho,R}(r)\frac{\mathcal{M}(r)}{g^{2}(r)}\leq 0
\end{aligned}
\end{equation}

Now we transplant the model space
solutions $\psi_{\rho,R}(r)$ of equation
(\ref{eqDir}) into the extrinsic annulus
$A_{\rho,R}=D_{R}(o)\setminus \bar D_{\rho}(o)$ in 
$S$ by defining
\[
\Psi_{\rho,R}\colon A_{\rho,R} \to \erre, \quad
\Psi_{\rho,R}(x)=\psi_{\rho,R}(r(x)).
\]
Here the extrinsic ball $D_{\rho}(o)$
is as in Definition (\ref{ExtBall}) and $D_{R}(o)$ is that 
component of $B_{R}(o)\cap S$ which contains $D_{\rho}(o)$.

Then, the hyptohesis (\ref{eqKcomp1}) on the sectional curvatures, and the assumption (\ref{geConvexity}) on the $o$-radial convexity leads to the following estimate using Proposition  \ref{corLapComp} (i), (recall that $\psi_{\rho,R}'(r) \geq 0$)
\begin{equation}
\begin{split}
\label{eqLap1}
\Delta^{S}\psi_{\rho,R}(r)  &\leq  \left(\psi_{\rho,R}(r) ''(r)
- \psi_{\rho,R}(r) '(r)\eta_{w}(r) \right)
 \Vert \nabla^{S} r \Vert^{2}
\\&+ m\,\psi_{\rho,R}(r) '(r) \left(\eta_{w}(r) - h(r)\right).
\end{split}
\end{equation}
so, using the second inequality in (\ref{parenth}) and that  $\Vert \nabla^{S}(r)\Vert \geq
g(r)$, we have:

\begin{equation}\label{laplacianIneq1}
\begin{aligned}
\Delta^{S}\psi_{\rho,R}(r(x)) \, &\leq \,
\left(\psi_{\rho,R}''(r(x))- \psi_{\rho,R}'(r(x))\eta_{w}(r(x))\right)\,g^{2}(r(x)) \\
& \phantom{mm}+ m\psi_{\rho,R}'(r(x))\left( \eta_{w}(r(x))
- h(r(x))\right)\, \\ &= \,g^{2}(r(x)) \Lmod_g \psi_{\rho,R}(r(x)) \, \\
&= \, 0 \, \\ &= \, \Delta^{S}v(x) \quad ,
\end{aligned}
\end{equation}
where $v(x)$ is the Laplace potential function
for the extrinsic annulus  $A_{\rho, R}\, = \,
D_{R} - D_{\rho}$, setting $v\vert_{\partial D_{\rho}}\,=\,0$ and
$v\vert_{\partial D_{R}}\,=\,1$.\\

Now, we apply the maximum principle to inequality (\ref{laplacianIneq1}) to obtain:
\begin{equation}
\psi_{\rho,R}(r(x)) \, \geq \, v(x)\,\, , \,\, \textrm{for
all} \, \, x \in A_{\rho, R} \quad .
\end{equation}

This implies in particular that on $\partial
D_{\rho}$ we have
\begin{equation}
\Vert \nabla^S\psi_{\rho,R}\Vert \, \geq \, \Vert
\nabla^{S}v(x)_{|_{\partial D_{\rho}}} \Vert
\end{equation}

Thern, using equation (\ref{CapEq}), we get
    \begin{equation}\label{eqCapac2}
    \begin{aligned} \C(A_{\rho, R}) &= \, \int_{\partial D_{\rho}}\Vert \nabla^{S}v(x) \Vert \, d\nu \\ &\leq
    \int_{\partial D_{\rho}}\Vert \nabla^S \Psi_{\rho,R}\Vert\, d\mu\\
&= \psi'_{\rho,R}(\rho)
\int_{\partial D_{\rho}}\Vert \nabla^S r\Vert\, d\mu \\
& = \frac{\C_{\Lmod_g}(A_{\rho, R}^{w})}
{\Vol(\partial D_{\rho}^{w})}
\int_{\partial D_{\rho}}\Vert \nabla^S r\Vert\, d\mu.
    \end{aligned}
    \end{equation}

On the other hand  $D_{\rho}(o)$ is precompact with a
smooth boundary and thence, 
\begin{equation}\label{PosInt}
\int_{\partial D_{\rho}}\Vert \nabla^S
r\Vert\, d\mu \, > 0 .
\end{equation}

Now, we have, using equations (\ref{eqModelCap}) and  (\ref{eqParabCond}):

\begin{equation}\label{intcond}
\begin{aligned}
&\C\bigl(\bar D_{\rho}(o), S^{m}\bigr)=\\ &
\lim_{R\to \infty}\C\bigl(\bar
D_{\rho}(o),D_{R}(o)\bigr)\,\\& \leq (\int_{\partial D_{\rho}}\Vert \nabla^S
r\Vert \, d\mu)\left(\lim_{R \to \infty}\frac{\C_{\Lmod_g}(A_{\rho, R}^{w})}
{\Vol(\partial D_{\rho}^{w})}\right)=0
\end{aligned}
\end{equation}

Thus, $D_{\rho}(o)$ is a precompact subset with
zero capacity in $S^{m}$, so the submanifold is parabolic.
\medskip

{\bf Assertion (A) in Theorem \ref{thmMain2}}.
\medskip
\\
Concerning the proof of assertion (A) in Theorem \ref{thmMain2} and under balance condition (\ref{balance2}), we have that the solution of the problem (\ref{eqDir}) satisifies:

\begin{equation}\label{parenth2}
\begin{aligned}
\psi'_{\rho,R}(r)& \geq 0\\
\psi''_{\rho,R}(r)& - \psi'_{\rho,R}(r)\eta_{w}(r) =- \psi'_{\rho,R}(r)\frac{\mathcal{M}(r)}{g^{2}(r)}\geq 0
\end{aligned}
\end{equation}

Then having into account that $\Vert \nabla^S r\Vert^2\geq g$ and $\phi_{\rho,R}''(r) - \phi'_{\rho,R}(r)\eta_{w}(r) \geq  0$ we obtain, applying Proposition \ref{corLapComp} (ii) to the transplanted solution $\psi_{\rho,R}$,
\begin{equation}\label{laplacianIneq2}
\begin{aligned}
\Delta^{S}\psi_{\rho,R}(r(x)) \, &\geq \,
\left(\psi_{\rho,R}''(r(x))- \psi_{\rho,R}'(r(x))\eta_{w}(r(x))\right)\,g^{2}(r(x)) \\
& \phantom{mm}+ m\psi_{\rho,R}'(r(x))\left( \eta_{w}(r(x))
- h(r(x))\right)\, \\ &= \,g^{2}(r(x)) \Lmod_g \psi_{\rho,R}(r(x)) \, \\
&= \, 0 \, \\ &= \, \Delta^{S}v(x) \quad ,
\end{aligned}
\end{equation}

We consider now inequality  (\ref{laplacianIneq2}), 
and procceed as in the proof of Theorem \ref{thmMain1}, but inverting the inequalities:
applying the maximum principle, we have on $\partial D_{\rho}$

\begin{equation*}
\Vert \nabla^S\psi_{\rho,R}\Vert \, \leq \, \Vert
\nabla^{S}v(x)_{|_{\partial D_{\rho}}} \Vert
\end{equation*}

 and using equation (\ref{CapEq}), we get
    \begin{equation*}
    %\begin{aligned} 
    \C(A_{\rho, R}) \geq
     \frac{\C_{\Lmod_g}(A_{\rho, R}^{w})}
{\Vol(\partial D_{\rho}^{w})}
\int_{\partial D_{\rho}}\Vert \nabla^S r\Vert\, d\mu.
   % \end{aligned}
    \end{equation*}

 Finally, having into account that we have inequality (\ref{PosInt}), and that condition $\int_{\rho}^{\infty} \Lambda_g(t)
\,  dt <  \infty$ is satisfied, we obtain 

\begin{equation}\label{intcond1}
%\begin{aligned}
\C\bigl(\bar D_{\rho}(o), S^{m}\bigr)\geq (\int_{\partial D_{\rho}}\Vert \nabla^S
r\Vert\, d\mu)\lim_{R \to \infty}\frac{\C_{\Lmod_g}(A_{\rho, R}^{w})}
{\Vol(\partial D_{\rho}^{w})}>0
%\end{aligned}
\end{equation}

Thus, $\bar D_{\rho}(o)$ is a compact subset with
positive capacity in $S^{m}$, and $S$ is hyperbolic.

%%%%%%%%%%%%%%%%%%%%%%%%%%%%%%%%%%%%%%%%%%%%%%%%%%%
%       Section 5: Proof of Theorem \ref{thmMain1} (ii)} and Theorem \ref{thmMain2} (ii)}
%%%%%%%%%%%%%%%%%%%%%%%%%%%%%%%%%%%%%%%%%%%%%%%%%%%
\section{Proof of assertion (B) in  Theorem \ref{thmMain1} and Theorem \ref{thmMain2}} \label{proofthmMain2}

Define now the following second order differential operator $\Lmod$,

\begin{equation}\label{operator}
\Lmod \phi(r) = \phi''(r) + \phi'(r)\left(
\mathcal{M}(r) -
\eta_{w}(r)\right).
\end{equation}

which is the same than in the proof above, with $g(r)=1 \,\,\forall r$, and consider the smooth radial solution $\phi_{\rho,R}(r)$ of the Dirichlet-Poisson problem associated to $\Lmod$ and defined on the interval $[\rho,R]$. As in the above proof, (see Proposition \ref{propDirSol}), it is easy to check that
\begin{equation} \label{eqPhi}
\phi_{\rho,R}(r) =  \frac{\int_{\rho}^{r}
\Lambda(t)\,dt}{\int_{\rho}^{R}\Lambda(t)\,dt}.
\end{equation}
where $\Lambda(r)$ is the function introduced in Definition \ref{Lambda} (i).

The corresponding 'drifted' $2$-capacity is
\begin{equation} \label{eqModelCap2}
\begin{aligned}
\C_{\Lmod}(A_{\rho, R}^{w})
&=\int_{\partial D_{\rho}^{w}}\langle\nabla^{M}\phi_{\rho,R},\nu\rangle\,dA\\
&=\Vol(\partial D_{\rho}^{w})\Lambda(\rho)\left(\int_{\rho}^{R}
\Lambda(t)\,dt\right)^{-1}.
\end{aligned}
\end{equation}
\medskip
{\bf Assertion (B) in Theorem \ref{thmMain1}}.
\\
Concerning assertion (B) in Theorem \ref{thmMain1}, we use equation (\ref{eqPhi}) and the balance condition (\ref{balance2}) to get
\begin{equation}\label{parenth2}
\begin{aligned}
\phi'_{\rho,R}(r)& \geq 0\\
\phi''_{\rho,R}(r)& - \phi'_{\rho,R}(r)\eta_{w}(r)\\&=- \phi'_{\rho,R}(r)\mathcal{M}(r) \geq 0
\end{aligned}
\end{equation}
because $\mathcal{M}(r)\leq 0$.
 
Having into account that $\Vert \nabla^S r\Vert^2\leq 1$ and $\phi_{\rho,R}''(r) - \phi'_{\rho,R}(r)\eta_{w}(r) \geq  0$ we obtain, applying Proposition \ref{corLapComp} (i)  to the transplanted function $\phi_{\rho,R}$,
\begin{equation}
\begin{aligned}
\Delta^{S}\phi_{\rho,R}(r(x)) \, &\leq \,
\left(\phi_{\rho,R}''(r(x))- \phi_{\rho,R}'(r(x))\eta_{w}(r(x))\right) \\
& \phantom{mm}+ m\phi_{\rho,R}'(r(x))\left( \eta_{w}(r(x))
- h(r(x))\right)\, \\ &= \Lmod\phi_{\rho,R}(r(x)) \, \\
&= \, 0 \, \\ &= \, \Delta^{P}v(x) \quad ,
\end{aligned}
\end{equation}
where $v(x)$ is the Laplace potential function
for the extrinsic annulus  $A_{\rho, R}\, = \,
D_{R} - D_{\rho}$, setting $v\vert_{\partial D_{\rho}}\,=\,0$ and
$v\vert_{\partial D_{R}}\,=\,1$.\\
Now parabolicity of $S$ follows as the proof of assertion (A) in Theorem \ref{thmMain1}.
\medskip
{\bf Assertion (B) in Theorem \ref{thmMain2}}.
\medskip
\\
To show assertion (B) in Theorem \ref{thmMain2}, we consider the same second order differential operator $\Lmod$, with the same smooth solution $\phi_{\rho,R}(r)$ to the same Dirichlet-Poisson problem defined on the interval $[\rho,R]$.

But now we have
\begin{equation}\label{parenth4}
\begin{aligned}
\phi'_{\rho,R}(r)& \geq 0\\
\phi''_{\rho,R}(r)& - \phi'_{\rho,R}(r)\eta_{w}(r) =- \phi'_{\rho,R}(r)\mathcal{M}(r)\leq 0
\end{aligned}
\end{equation}
because $\mathcal{M}(r) \geq 0$.

Then having into account that $\Vert \nabla^S r\Vert^2\leq 1$ and $\phi_{\rho,R}''(r) - \phi'_{\rho,R}(r)\eta_{w}(r) \leq  0$ we obtain, applying Proposition \ref{corLapComp} (ii) to the tranplanted function $\phi_{\rho,R}$, 
\begin{equation}
\begin{aligned}
\Delta^{S}\phi_{\rho,R}(r(x)) \, &\geq \,
\left(\phi_{\rho,R}''(r(x))- \phi_{\rho,R}'(r(x))\eta_{w}(r(x))\right) \\
& \phantom{mm}+ m\phi_{\rho,R}'(r(x))\left( \eta_{w}(r(x))
- h(r(x))\right)\, \\ &= \Lmod\phi_{\rho,R}(r(x)) \, \\
&= \, 0 \, \\ &= \, \Delta^{S}v(x) \quad ,
\end{aligned}
\end{equation}

Now hyperbolicity of $S$ follows as the proof of assertion (A) of Theorem \ref{thmMain2}.

%\begin{acknowledgements}\label{ackref}
%We would like to acknowledge professor Steen
%Markvorsen their useful comments concerning these results.
%\end{acknowledgements}
%\bibliographystyle{amsalpha}

%\end{document}

\enddocument